\documentclass[technote]{IEEEtran}
\usepackage{graphicx}
\usepackage{cite}
\usepackage{amsmath}
\usepackage{amsthm}
\usepackage{amssymb}
\usepackage{array}
\usepackage[caption=false]{subfig}
\usepackage[T1]{fontenc}
\usepackage[utf8]{inputenc}

\theoremstyle{definition}
\newtheorem{definition}{Definition}[section]

\begin{document}
\title{Chaotic Logistic Map Forecast using Fuzzy Time Series}

\author{\IEEEauthorblockN{Lucas Vinícius Ribeiro Alves}
\IEEEauthorblockA{Technical College\\
Federal University of Minas Gerais\\
E-mail: lucasvra@ufmg.br}
}

\maketitle

\begin{abstract}
This paper deals with the problem of forecast the Logistic Chaotic Map using Fuzzy Times Series (FTS). Chaotic Systems are very sensible to changes in its parameters and in the initial conditions, turning them into hard systems to model and forecast. In this case, we relay in the robustness of Fuzzy Time Series to model and forecast the logistic map.
We use the Akaike Information Criterion (AIC) as an index to determine the number of sub intervals for the definition of the fuzzy set.
\end{abstract}


\IEEEpeerreviewmaketitle

%
\section{Introduction}

The forecast of chaotic dynamics is a challenging task given the particularities of chaotic systems, like non periodic behavior and sensibility to initial conditions \cite{Lisi2001}. In this paper we deal with the problem of few steps forecast for the logistic map in chaotic behavior, using Fuzzy Time Series \cite{Song1993, Song1993U, Song1994, Song1997, chen2004, Chen2006, iqbal2020}. Chaotic time series prediction can be used as a benchmark for forecast methods, indication how they perform when applied to real data, like stock market data, weather forecast and bankruptcy prediction. 

Here we expect the robustness of the fuzzy time series models improve the performance of the forecast process, given the unpredictability of chaotic systems. The main objective of this paper is to determine if a first order fuzzy time series is capable of model the dynamics of a nonlinear dynamics with chaotic behavior. 

The paper is organized as follows: Section II presents the background information about the logistic map and fuzzy time series; Section III explains the methodology used in this paper and the procedures used to create the fuzzy time series model. Section IV presents the experimental results and a discussion about each experiment. The conclusions and directions for future work are presented in Section V.

%
\section{Background}

\subsection{Logistic Map and Chaotic Systems}

The logistic map is, usually, used to model population growth when the environment applies a saturation to the growth \cite{Kunsch2006}. The logistic map is defined by the following recurrence relation:

$$ x_{k+1} = r x_k(1-x_k)$$

Where the parameter $r$ is in the interval $[0,4]$. The logistic map is a population growth model for closed population living in a environment with limited resource for the subsistence of the population members.

\begin{figure}[htbp]
    \centering
    \includegraphics[width=0.9\columnwidth]{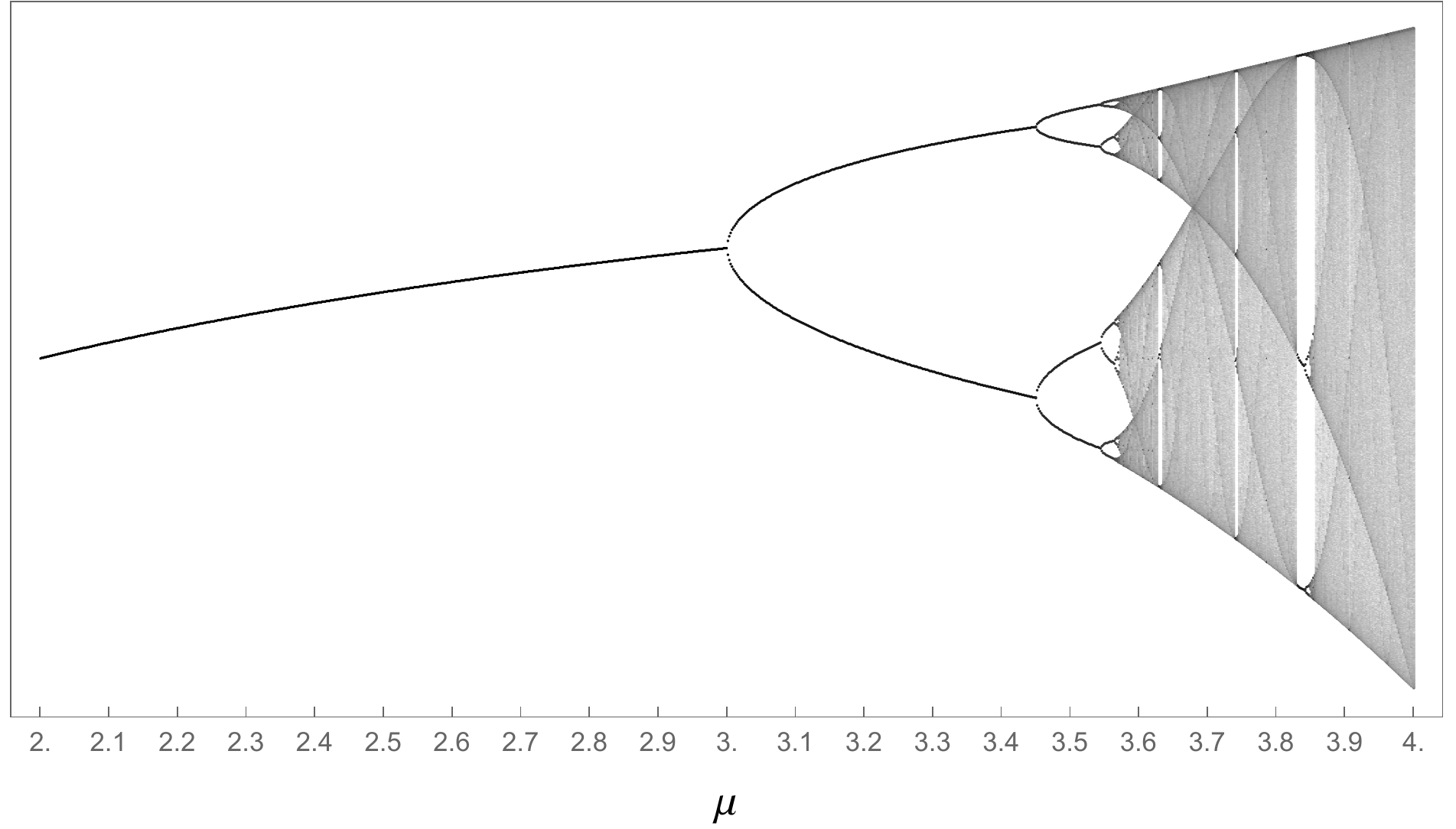}
    \caption{Bifurcation diagram of the logistic map}
    \label{fig:bifurc}
\end{figure}

Changing the parameter $r$ may imply in drastic changes in the dynamics of the system, in Fig.~\ref{fig:bifurc} is shown the bifurcation diagram of the logistic function. As we may see, for $r>3.6$ the system exhibits, mostly, chaotic behavior.

\begin{figure}[htbp]
    \centering
    \includegraphics[width=0.9\columnwidth]{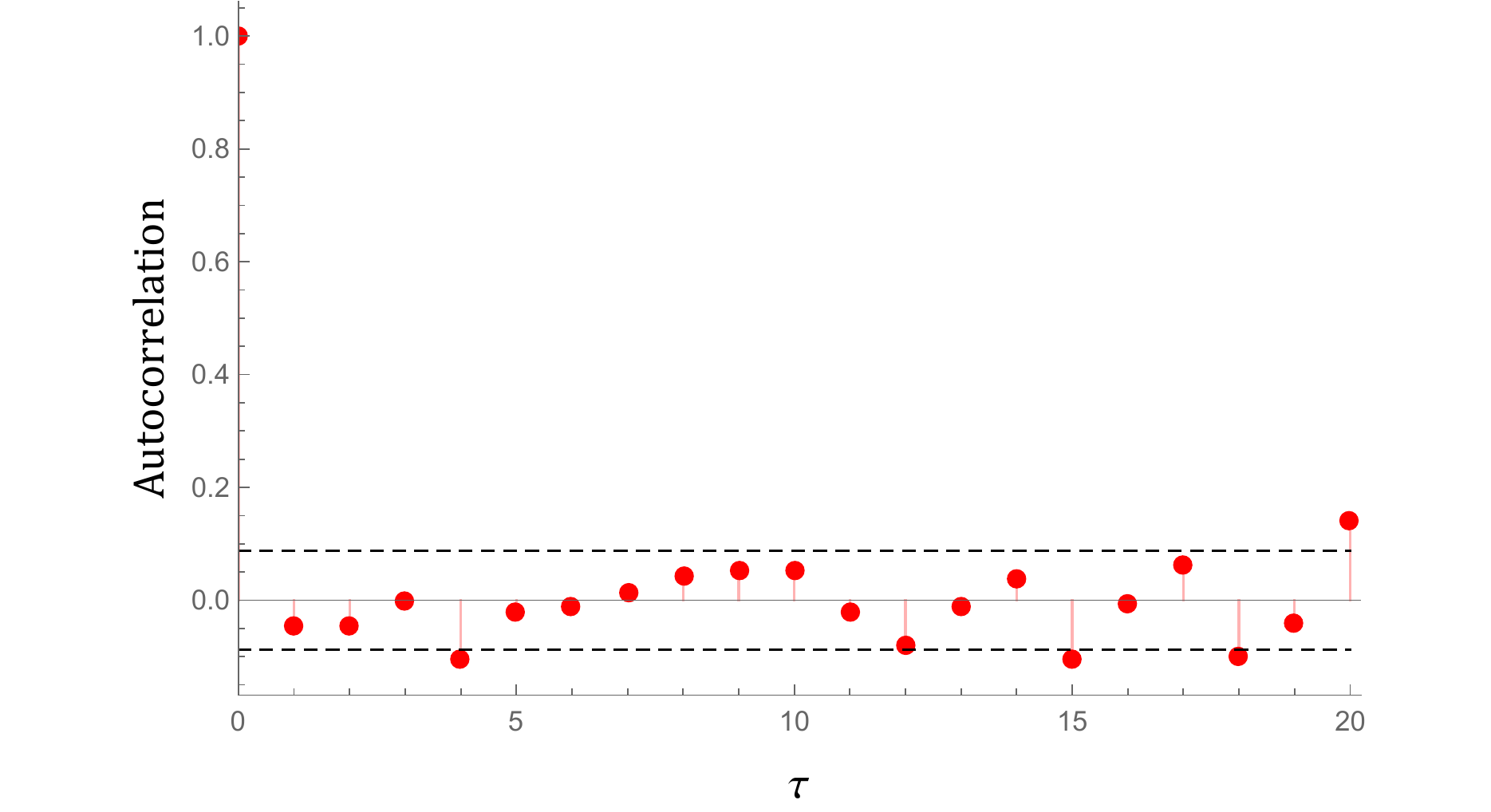}
    \caption{Autocorrelation function of the logistic map (removed the mean), with $x_1=0.2$ and $r=3.999$}
    \label{fig:acf}
\end{figure}

As we may see in the Fig.~\ref{fig:acf}, autocorrelation function of the logistic map, in chaotic behavior there is no statistically significant correlation between the time series and a delayed version. With no significant correlation, linear models are incapable of forecast the system well.

\begin{figure}[htbp]
    \centering
    \includegraphics[width=0.9\columnwidth]{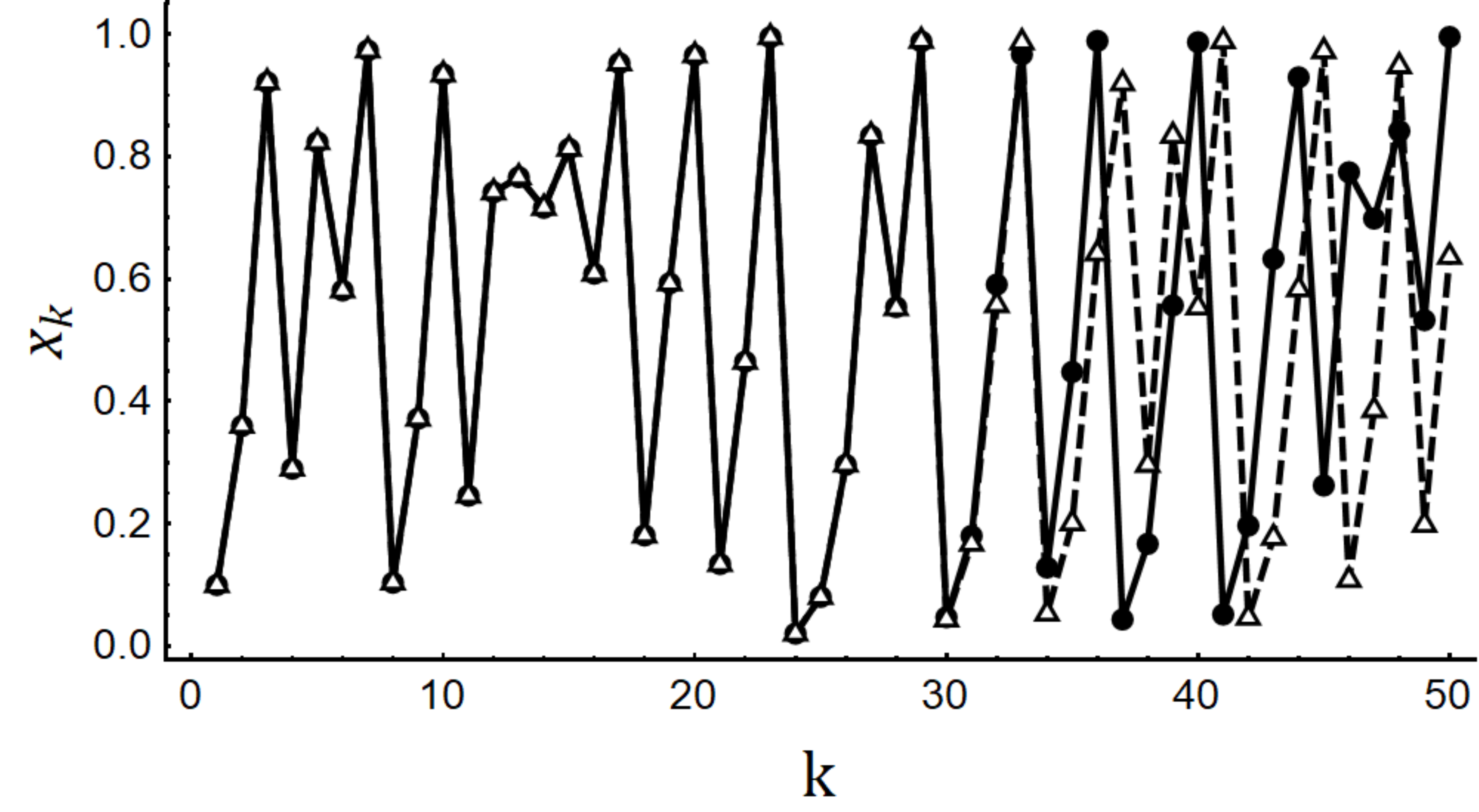}
    \caption{Two simulations of the logistic map, one with $x_1=0.1$ and $r=3.999$ and the other with $x_1=0.10000000001$ and $r=3.999$}
    \label{fig:sim}
\end{figure}

An example of the sensibility to the initial conditions is shown on Fig.~\ref{fig:sim}, where a difference of $10^-11$ in the initial conditions leads the system to different trajectories after about 35 iterations.

Other fuzzy approaches were already used to forecast chaotic systems, based on fuzzy neural networks \cite{Zhang2008, Yang2006}.

\subsection{Fuzzy Time Series}

The concept of Fuzzy Times Series was developed by \cite{Song1993} in order to model a dynamical process where the observations are linguistic values. The main characteristics of Fuzzy Time Series (FTS) are:

\begin{itemize}
    \item FTS are dynamical processes.
    \item The observations of FTS are fuzzy sets.
    \item The universe o discourse for the fuzzy sets are subsets of $\mathbb{R}^1$.
    \item Conventional time series models are not applicable to these processes.
\end{itemize}

A formal definition of fuzzy time series is shown below:

\begin{definition}
Let $Y(t)$, $(t=\ldots,1,2,3,\ldots)$, a subset of $\mathbb{R}^1$, be the universe of discourse on which fuzzy sets $f_i(t)$ $(i = 1,2,3, \ldots)$ are defined and $F(t)$ is the collection of $f_i(t)$ $(i = 1,2,3, \ldots)$. $F(t)$ is then called a Fuzzy Time Series on $Y(t)$, $(t=\ldots,1,2,3,\ldots)$.
\end{definition}

%
\section{Methodology}

In this section we first present a convectional method for forecast using fuzzy time series, then we show the proposed use of the Akaike Information Criterion as a tool to determine the number of intervals used to create the fuzzy sets.

%
\subsection{Convectional Fuzzy Time Series Method}
The most basic method to forecast using Fuzzy Time Series was present by Song and Chissom in \cite{Song1993} where they forecast the enrollments of a university \cite{Song1993U,Song1994}. First, we must have some definitions \cite{Song1993}:

\begin{definition} 
If there is a fuzzy relationship $R(t, t-1)$, such that $F(t)=F(t-1)\times R(t,t-1)$, where $\times$ is an operator, the $F(t)$ is caused by $F(t-1)$. The relationship between $F(t)$ and $F(t-1)$ is denoted by $F(t-1) \to F(t)$.
\end{definition}

\begin{definition} 
Denoting $F(t-1)$ by $A_i$ and $F(t)$ by $A_j$, the relationship between $F(t-1)$ and $F(t)$ can be defined as the logical relationship $A_i \to A_j$.
\end{definition}

\begin{definition} 
Fuzzy logical relationships with the same left-hand sides can be grouped into fuzzy logical relationship groups. So, $\{A_i \to A_1, A_i \to A_2, \ldots \}$ can be denoted as $A_i \to A_1, A2, \ldots$.
\end{definition}

\noindent
Then, the method can be summarized in 8 steps:

\begin{itemize}
    \item[Step 1] Define universe of discourse (universal set $U$) using the minimum a maximum on the time series
    \item[Step 2] Partition $U$ into intervals with the same length.
    \item[Step 3] Define fuzzy sets $A_i$ using the intervals.
    \item[Step 4] Fuzzify the data.
    \item[Step 5] Determine the fuzzy logical relationships $A_i \to A_j$.
    \item[Step 6] Group the fuzzy logical relationships with the same left-hand side and calculate the matrix $R_i$ for each fuzzy relationship group.
    \item[Step 7] Forecast.
    \item[Step 8] Defuzzify to obtain the desired result.
\end{itemize}

%
\subsection{Interval Definition}

The definition of the intervals used in the fuzzification process is a crucial step in the fuzzy time series modeling process \cite{Huarng2001}. The length of the intervals significantly affects the forecast result, becoming an area of study in the field of fuzzy time series. There are many approaches to determine an effective interval length, as using the distribution or average of the differences \cite{Huarng2001}, using optimization \cite{Egrioglu2010}, using clustering \cite{Zhang2012}

\subsubsection{Interval Definition using Akaike Information Criterion}

In this work, the observations of the logistic map are confined in the interval $[0,1]$, so, in order to determine the intervals, we may only choose the number of intervals. To do so, we obtain a model for many numbers of interval and apply the Akaike Information Criterion (AIC) \cite{Brockwell2002}, using the residuals. The number of intervals, or the length of the intervals, is a parameter of the system, and should be used with parsimony. Usually, increasing the number of intervals reduce the forecast error, but many intervals may lead to over fitting.

To calculate the Akaike Information Criterion we use the following formula:

$$ AIC(n) = N ln(RSS/N) + 2 n $$

\begin{figure}
    \centering
    \includegraphics[width=0.8\columnwidth]{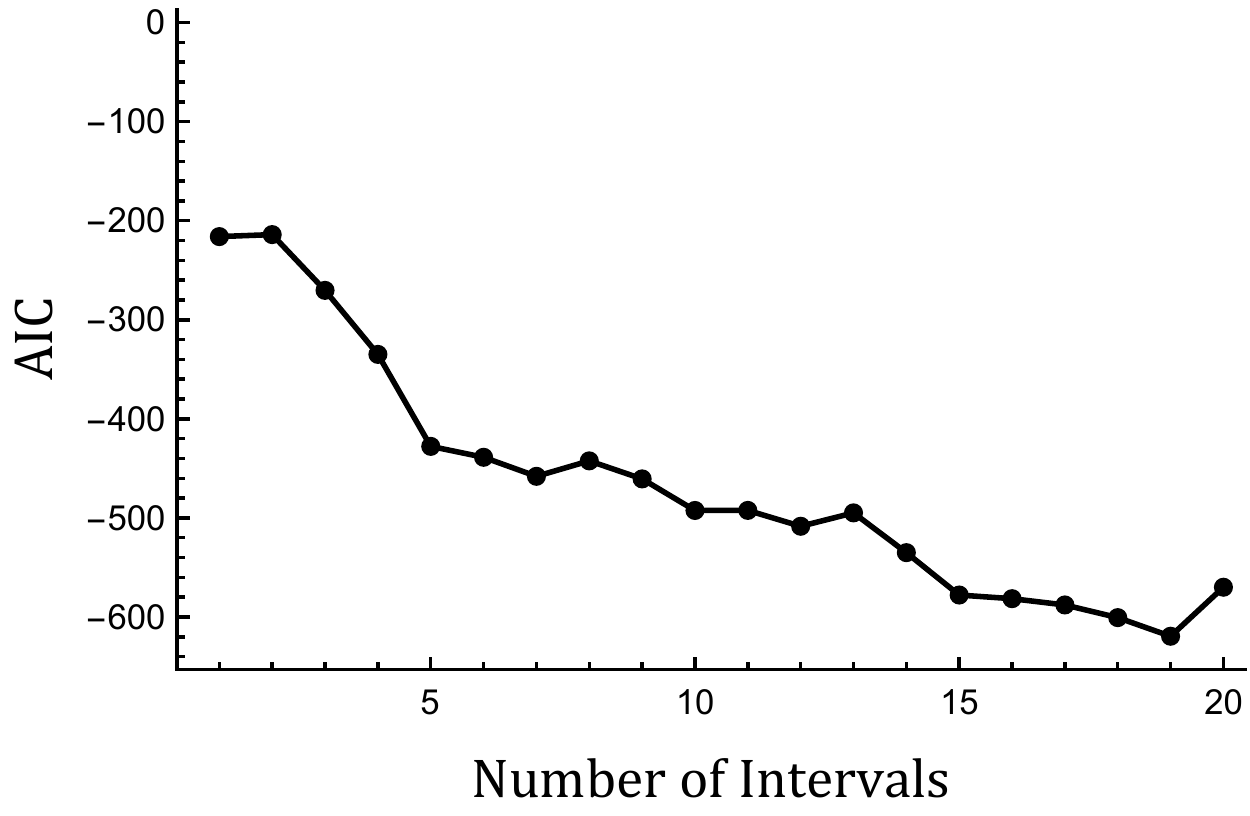}
    \caption{Plot of the Akaike Information Criterion for the logistic map, with $r = 3.999$ and $x_1 = 0.1$}
    \label{fig:aic}
\end{figure}

Where $n$ is the number of sub intervals, $N$ is the number of samples and $RSS$ is the residual sum of squares. The number of intervals used is the first local minimum of the AIC plot. Fig.~\ref{fig:aic} shows the AIC plot for the logistic map, with $r = 3.999$, $x_1 = 0.1$, and a sample size $N = 100$. As we can see, the first minimum happens in $n=7$ intervals. Increasing the number of intervals improves the model but may lead to over fitting and increases the model complexity.

\subsubsection{Interval Definition using the Average Based Method}

The average based method is based on the average of the first differences of data on the time series \cite{Huarng2001}. The algorithm for the definition of the length is:

\begin{itemize}
    \item[Step 1] Calculate the absolute differences between $x_{k+1}$ and $x_k$ and take the average of the resulting series.
    \item[Step 2] Take one half of the average obtained as the length of the interval.
    \item[Step 3] Using the length obtained, determine a base for the length. Let $\mu$ be the average, if $10^n < \mu \leq 10^{n+1}$, then $10^n$ is the base.
    \item[Step 4] Given the base, round the length according to the base.
\end{itemize}

For the logistic map, with $r = 3.999$, $x_1 = 0.1$, and a sample size $N = 100$, half of the average of the absolute differences is $0.210145$, so the base is $0.1$, leading to a length of $0.2$ and 5 intervals.

%
\subsection{Algorithm Application}
\begin{figure}
    \centering
    \includegraphics[width=0.8\columnwidth]{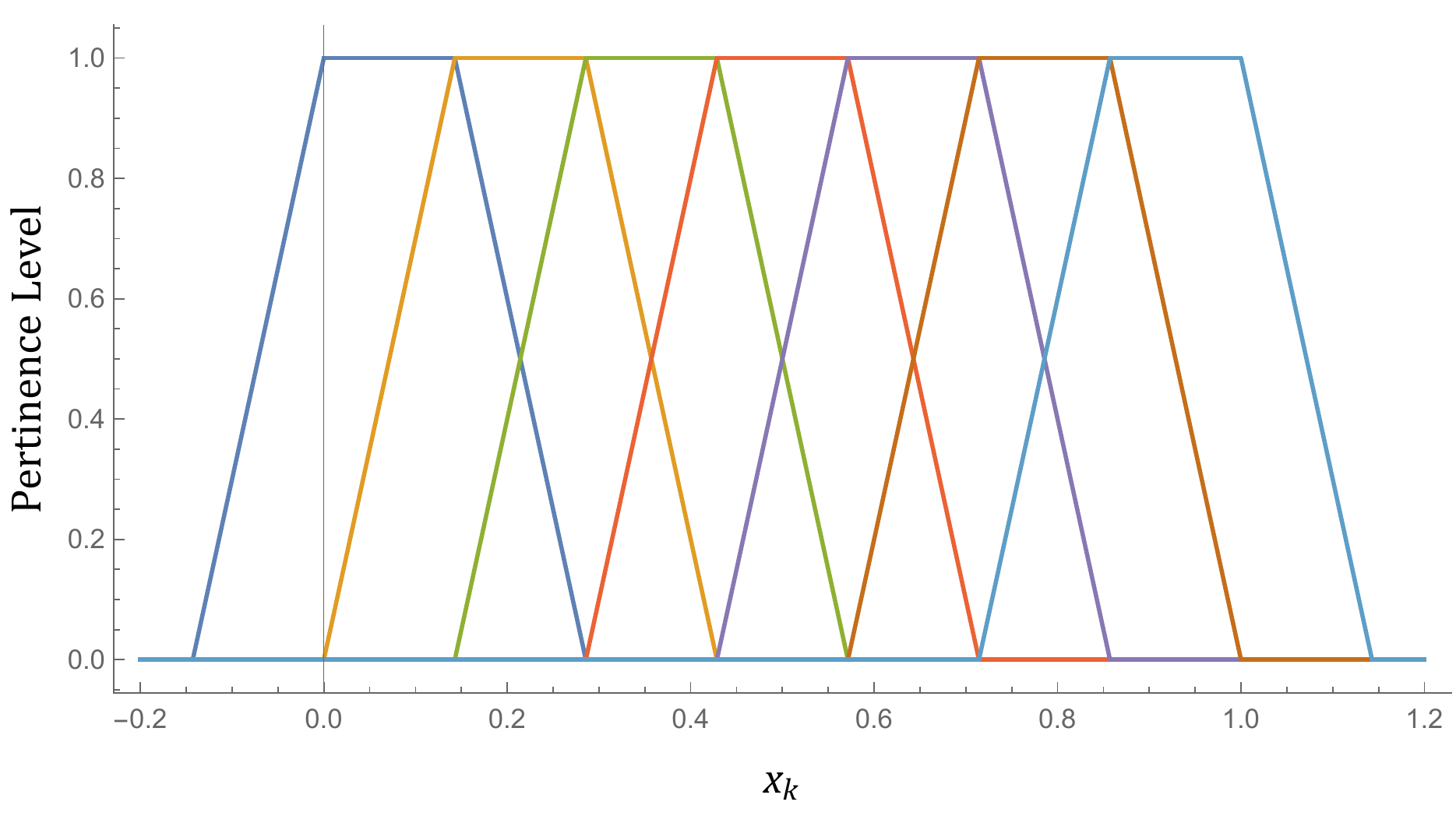}
    \caption{Example of set of fuzzy membership functions}
    \label{fig:pert}
\end{figure}

The first step is to generate the fuzzy membership function to each interval $[a,b]$ as a trapezoidal function, as shown on Fig.~\ref{fig:pert}. The trapezoidal membership function is defined by the following formula:

$$\mu (x) = 
\begin{cases} 
0                    & \text{if } x \leq 2a-b \\
\frac{x-2a+b}{b-a}   & \text{if } 2a-b \leq x \leq a \\
1                    & \text{if } a \leq x \leq b \\
\frac{3b-3a-x}{b-a}  & \text{if } b \leq x \leq 2b-a \\
0                    & \text{if } 2b-a \leq x
\end{cases}
$$

After that, we evaluate the observations of the system in each fuzzy set, generating a table of pertinence to each fuzzy set in each time. On TABLE~\ref{tab:pert} is shown the pertinence on each fuzzy set for the first 10 samples of the time series.

Using the table, we determine the fuzzy relationships using taking the fuzzy sets with pertinence 1 in subsequent time. So, if in sample $k$ the fuzzy set $A_i$ has pertinence 1 and in sample $k+1$ the fuzzy set $A_i$ has pertinence 1, we may define a fuzzy relation $A_i \to A_j$.

\begin{table}[htbp]
\centering
\caption{Pertinence on each fuzzy set for samples 1 to 10}
\footnotesize
\begin{tabular}{|c|c|c|c|c|c|c|c|c|c|}\hline
$k$  & $A_1$ & $A_2$ & $A_3$ & $A_4$ & $A_5$ & $A_6$ & $A_7$ \\\hline
1  & 0.6 & 1.0 & 0.4 & 0.0 & 0.0 & 0.0 & 0.0 \\ \hline
2  & 0.0 & 0.0 & 0.0 & 0.5 & 1.0 & 0.5 & 0.0 \\ \hline
3  & 0.0 & 0.0 & 0.0 & 0.0 & 0.0 & 0.6 & 1.0 \\ \hline
4  & 0.0 & 0.8 & 1.0 & 0.2 & 0.0 & 0.0 & 0.0 \\ \hline
5  & 0.0 & 0.0 & 0.0 & 0.0 & 0.0 & 0.9 & 1.0 \\ \hline
6  & 0.0 & 0.0 & 0.7 & 1.0 & 0.3 & 0.0 & 0.0 \\ \hline
7  & 0.0 & 0.0 & 0.0 & 0.0 & 0.0 & 0.0 & 1.0 \\ \hline
8  & 1.0 & 0.1 & 0.0 & 0.0 & 0.0 & 0.0 & 0.0 \\ \hline
9  & 1.0 & 0.2 & 0.0 & 0.0 & 0.0 & 0.0 & 0.0 \\ \hline
10 & 1.0 & 0.8 & 0.0 & 0.0 & 0.0 & 0.0 & 0.0 \\ \hline
\end{tabular}
\label{tab:pert}
\end{table}

After finding the fuzzy relationships for each consecutive pair of samples, we remove the duplicates and group then by the first member in the relationship. 

\begin{center}
\begin{tabular}{l}
 $A_1 \to$ \{$A_1$,$A_2$,$A_3$,$A_4$\} \\
 $A_2 \to$ \{$A_4$,$A_5$,$A_6$\} \\
 $A_3 \to$ \{$A_6$,$A_7$\} \\
 $A_4 \to$ \{$A_7$\} \\
 $A_5 \to$ \{$A_6$,$A_7$\} \\
 $A_6 \to$ \{$A_4$,$A_5$,$A_6$\} \\
 $A_7 \to$ \{$A_1$,$A_2$,$A_3$,$A_4$\} \\
\end{tabular}
\end{center}

We must then define the operator $\times$ between two vectors:

\begin{definition}
Let $C$ and $B$ be two row vectors of dimension $n$ and let $D=(d_{ij})=C^T \times B$, where $d_{ij} = min(c_i, b_j)$ where $c_i$ is the $i$-th element of $C$ and $b_j$ is the $j$-th element of $B$. 
\end{definition}

So, given the fuzzy logical relationships, for each relationship $A_i \to A_j$ we should compute $R_k = A_i^T \times A_j$. Here, the approach diverges from \citeonline{Song1993}, instead of using a default vector for $A_i$, we use the first vector in the time series with maximum pertinence in the fuzzy set $A_i$.

In the example above, one should compute $m = 19$ matrices $R_k$ in order to obtain:

$$ R = \bigcup\limits_{k=1}^{m} R_k $$

Using $R$, the forecast model is given by $A_{k+1} = A_k \circ R $, where $A_k$ is the fuzzified observation in time $k$ and $A_{k+1}$ is the fuzzified forecast. An example of matrix $R$ is shown bellow:

$$ R = \left(
\begin{array}{ccccccc}
 1.0 & 1.0 & 1.0 & 1.0 & 0.8 & 0.5 & 0.2 \\
 0.4 & 0.4 & 0.9 & 1.0 & 1.0 & 1.0 & 0.8 \\
 0.0 & 0.0 & 0.1 & 0.4 & 0.7 & 1.0 & 1.0 \\
 0.0 & 0.0 & 0.0 & 0.0 & 0.0 & 0.5 & 1.0 \\
 0.0 & 0.0 & 0.0 & 0.0 & 0.7 & 1.0 & 1.0 \\
 0.3 & 0.6 & 0.7 & 1.0 & 1.0 & 1.0 & 0.8 \\
 1.0 & 1.0 & 1.0 & 1.0 & 0.5 & 0.4 & 0.2 \\
\end{array}
\right)$$

If necessary, a defuzzification can be done in order to obtain the forecast in numeric value.

%
\subsection{Defuzzification}

There are many different approaches to defuzzify data in the fuzzy set theory, being a way to produce quantifiable value in crisp logic given determined fuzzy sets and pertinence of an observation to each set. 

In the field of fuzzy time series, the most convectional method, also used in this paper, is the following \cite{Song1993}:

\begin{enumerate}
    \item If the vector $A_k$ has only one maximum, the defuzzified forecast is the midpoint of the interval corresponding to the maximum.
    \item If the vector $A_k$ has two or more consecutive maximums, the defuzzified forecast is the midpoint of the conjunct intervals corresponding to the maximums.
    \item Otherwise, standardize the fuzzy forecast vector $A_k$ and use the midpoint of each interval to calculate the centroid of the fuzzy set as the defuzzified forecast.
\end{enumerate}

%
\section{Results and Discussion}

In the tests we are going to use four models:

\begin{itemize}
    \item[Model 1] A first order fuzzy time series model.
    \item[Model 2] A linear autoregressive model, in the form \\$x_{k+1} = \theta x_k$
    \item[Model 3] A quadratic autoregressive model, in the form \\$x_{k+1} = \theta x^2_k$
    \item[Model 4] A combination of models 2 and 3, in the form \\$x_{k+1} = \theta_1 x^2_k + \theta_2 x_k$
\end{itemize}

To exemplify each model, we used two time series, the first one is generated by the logistic map with $r=3.999$ and $x_1 = 0.1$ and is used to obtain the model and a second time series, generated by the logistic map with $r=3.999$ and $x_1 = 0.2$.

It is easy to see that model 4 is the exact logistic map, so it should give the better results. Model 2 is linear and will not generate good forecasts, since the logistic map has no significant autocorrelation with a delayed version of itself. 

\begin{figure}[htbp]
    \centering
    \subfloat[Model 1 (FTS)]{
       \includegraphics[width=0.8\columnwidth]{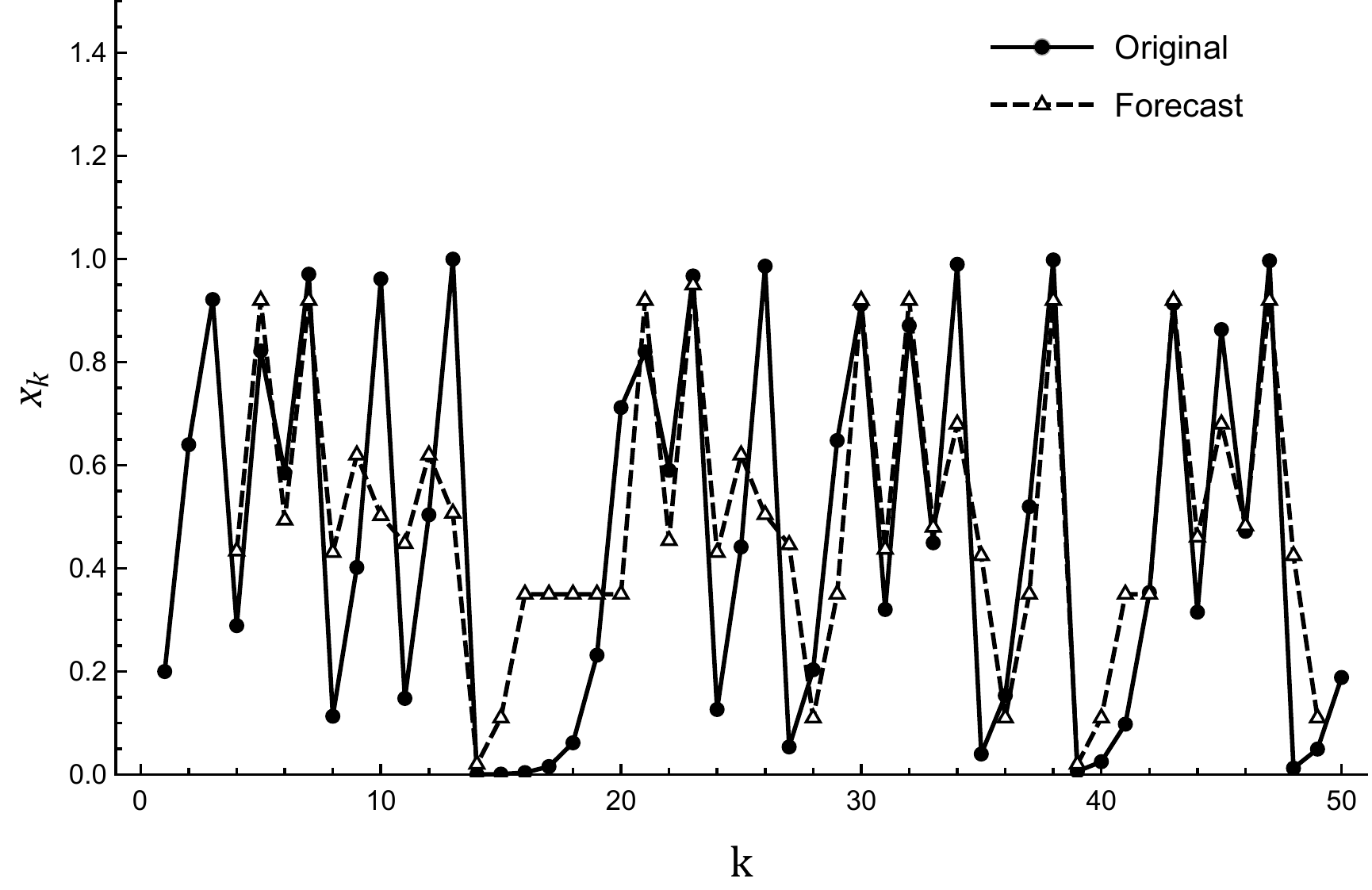}
     }
     
     \subfloat[Model 2]{
       \includegraphics[width=0.8\columnwidth]{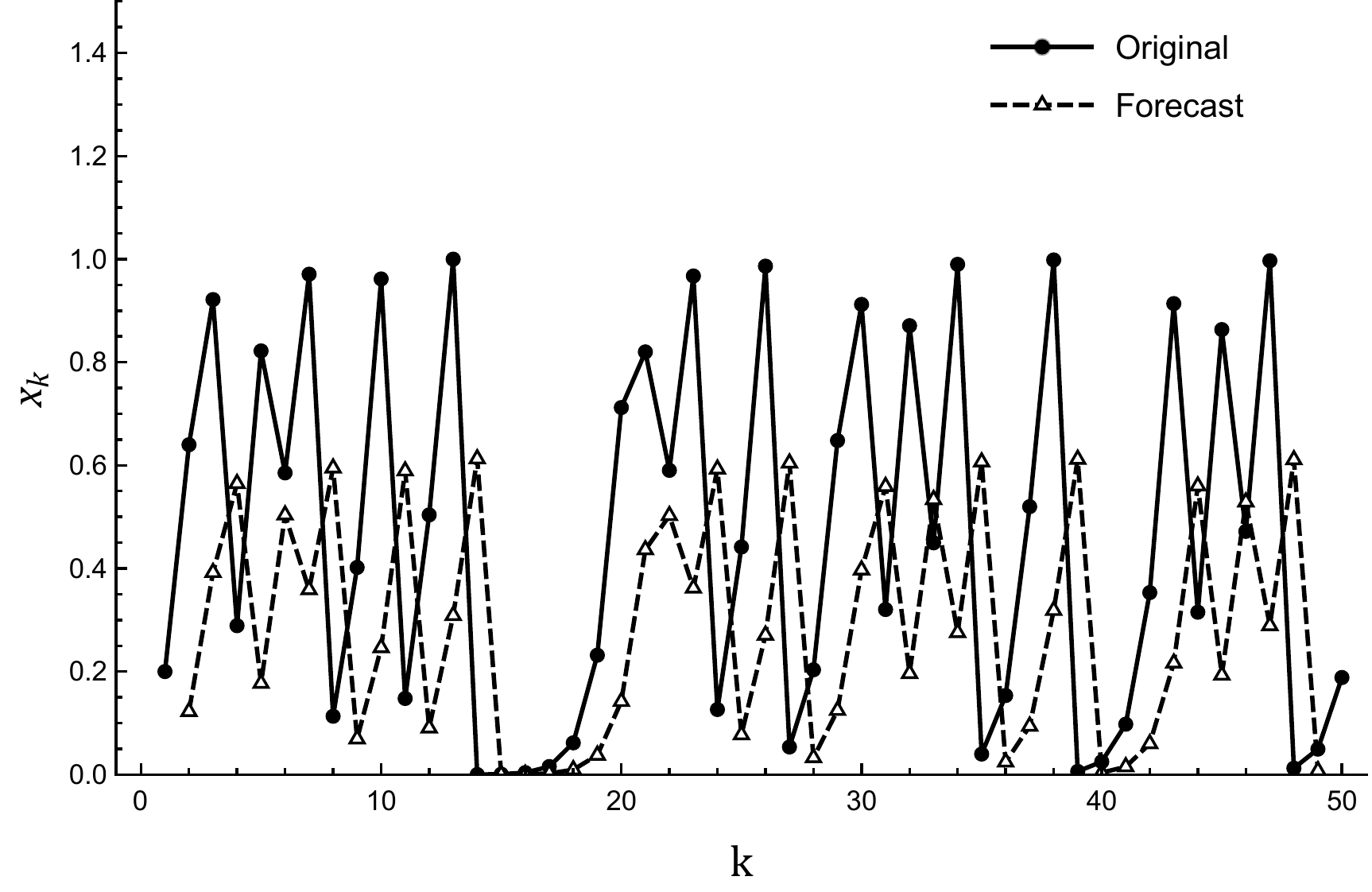}
     }
     
     \subfloat[Model 3]{
       \includegraphics[width=0.8\columnwidth]{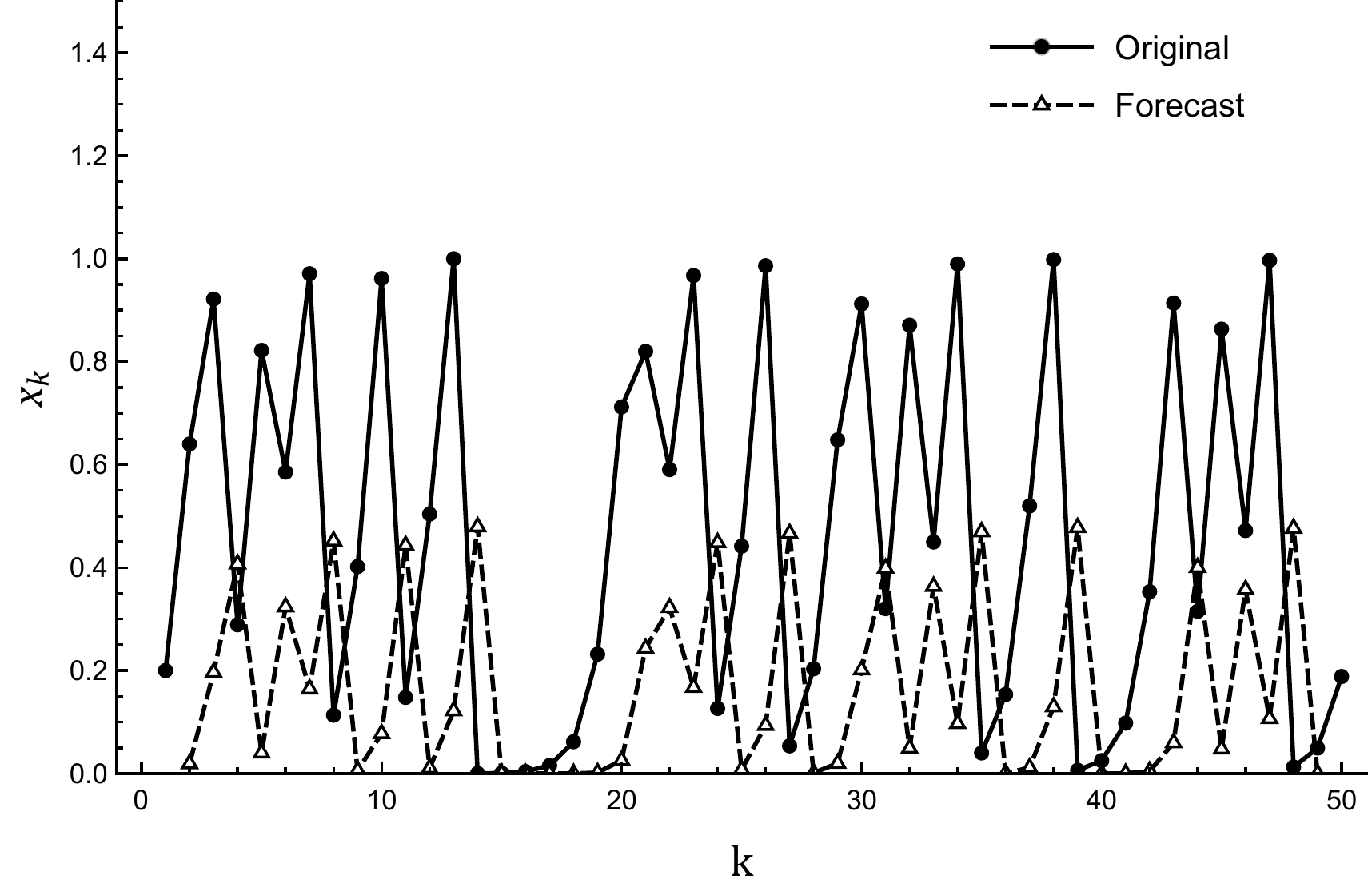}
     }
     
     \subfloat[Model 4]{
       \includegraphics[width=0.8\columnwidth]{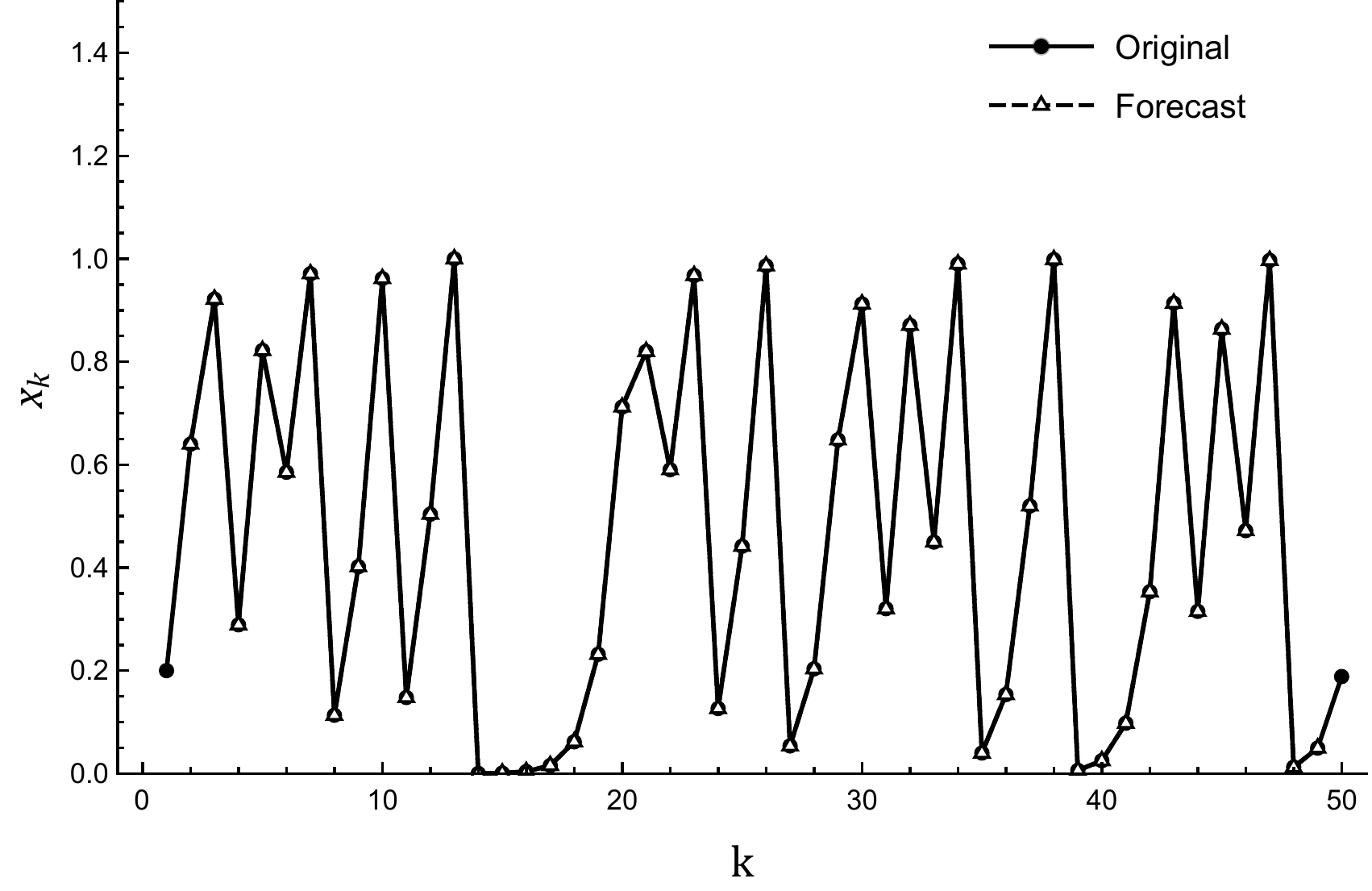}
     }
     \caption{One step ahead forecast using the four models}
     \label{fig:test1}
\end{figure}

On Fig.~\ref{fig:test1} is shown the one step ahead forecast of each model compared to the original time series. As can be seen, the model 4 perfectly follows the system, what was expected. The fuzzy time series model (model 1) is far better than model 2 and model 3 the forecast.

\begin{figure}[htbp]
    \centering
    \includegraphics[width=0.9\columnwidth]{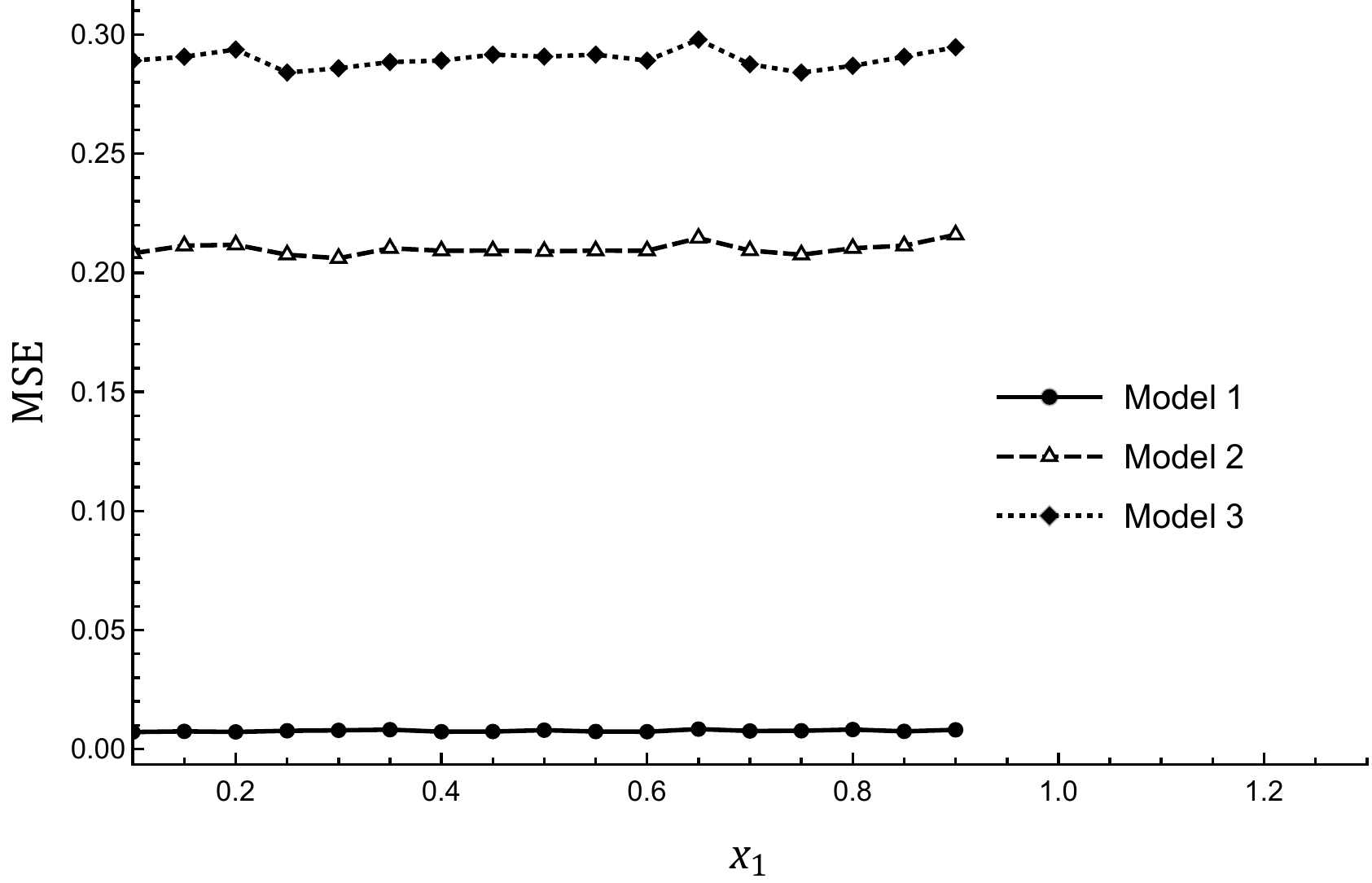}
    \caption{Forecast mean squared error of the models 1, 2 and 3 for different initial conditions ($x_1$)}
    \label{fig:res1}
\end{figure}

The first test applied was compare the MSE for one step ahead forecast of the logistic map, the time series used was generated using $r=3.999$ and $x_1 \in [0.1,0.9]$ with 1000 observations. The first 500 samples were used to estimate the model and the results shown on Fig.~\ref{fig:res1}, are the MSE (mean squared error) of the one step ahead forecast for the final 500 samples in each initial condition. The model 4 was omitted because its error is close to zero since it is the exact model.

As we can see, the MSE of the fuzzy time series model is much better than the other two models, in this case were used 9 intervals to define the fuzzy sets, for all initial conditions. On TABLE~\ref{tab:res1} is shown the mean and the variance of the forecast MSE when $x_1=0.2$, showing how the fuzzy time series model has a small mean and variance of the error.

\begin{table}[htbp]
    \centering
    \begin{tabular}{c|rr}
  & \text{Mean} & \text{Variance} \\\hline
 \text{Model 1} & \textbf{0.007237} & \textbf{0.0001187} \\
 \text{Model 2} & 0.211746 & 0.0292752 \\
 \text{Model 3} & 0.293699 & 0.0687509 \\
    \end{tabular}
    \caption{Forecast mean squared error of the models 1, 2 and 3 $x_1=0.2$}
    \label{tab:res1}
\end{table}

As a second test, we compare the MSE of the three models for the one step ahead forecast of the time series generated using $x_1 = 0.1$ and $r \in [3,4]$ with 1000 observations. Again, the first 500 observations are used to estimate the models and the other 500 observations are used to validation. Here, for each value of $r$ the AIC was used to determine the number of intervals used in the model 1. 

\begin{figure}[htbp]
    \centering
    \includegraphics[width=0.9\columnwidth]{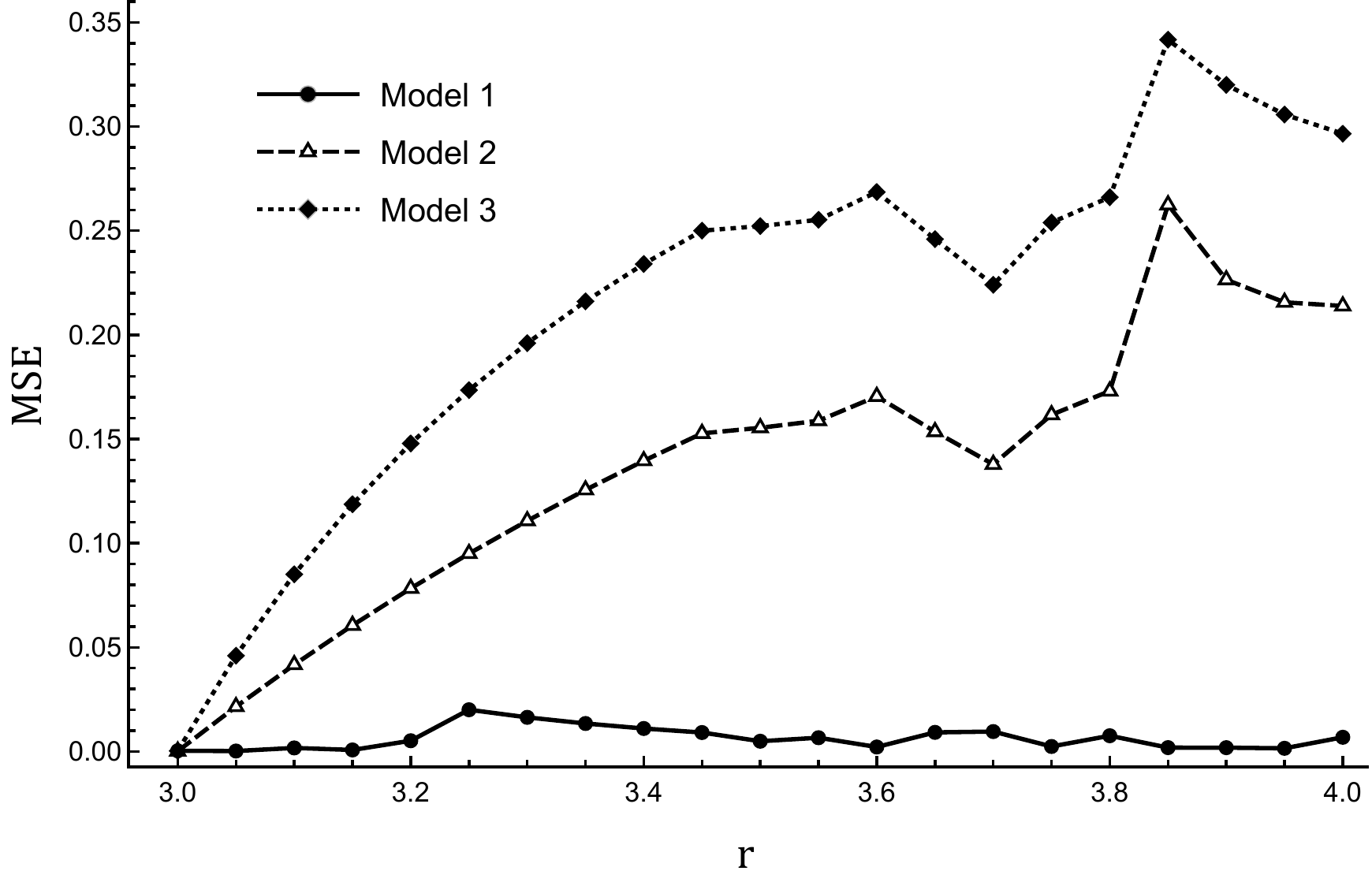}
    \caption{Forecast mean squared error of the models 1, 2 and 3 for different values of $r$}
    \label{fig:res2}
\end{figure}

As we may see on Fig.~\ref{fig:res2}, the fuzzy time series model performed better in the one step ahead forecast and is more insensitive to variations in the parameter $r$ than the other two models.

TABLE~\ref{tab:res2} summarizes the same test for $r=3.999$, with the information about the variance of the error. The fuzzy time series model was the model with the smallest error and variance.

\begin{table}[htbp]
\caption{Forecast mean squared error of the models 1, 2 and 3 for $r=3.999$}
    \centering
    \begin{tabular}{c|rr}
 & \text{Mean} & \text{Variance} \\\hline
 \text{Model 1} & \textbf{0.011129} & \textbf{0.0001643} \\
 \text{Model 2} & 0.208182 & 0.0284782 \\
 \text{Model 3} & 0.289097 & 0.0641057 \\
    \end{tabular}
    \label{tab:res2}
\end{table}

\begin{figure}[htbp]
    \centering
    \includegraphics[width=0.8\columnwidth]{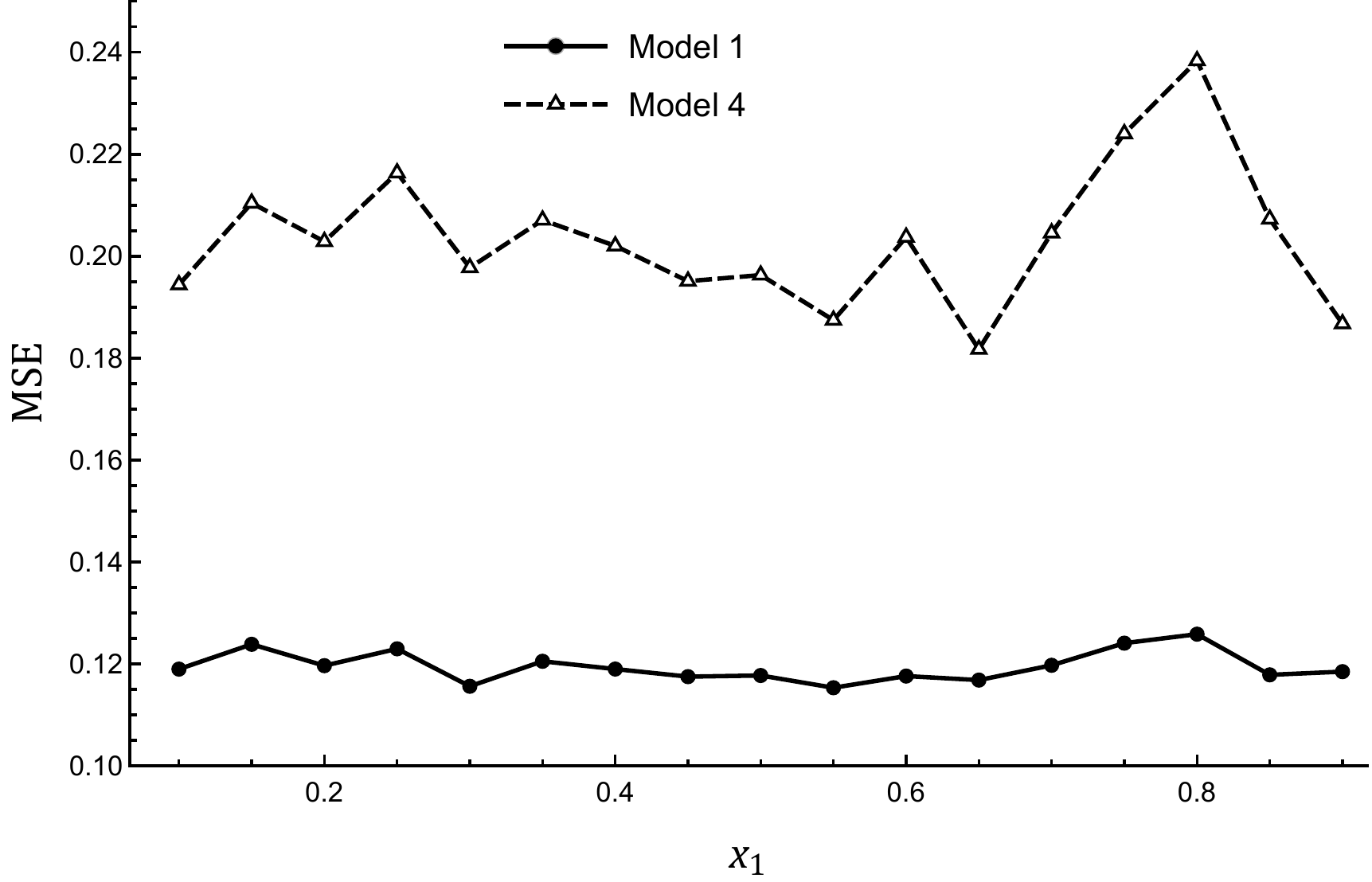}
    \caption{Forecast error of the models 1 and 4 for three steps ahead and different initial conditions ($x_1$)}
    \label{fig:res3}
\end{figure}

As a third test, we generate 1000 observations of the logistic map with $r=3.999$ and $x_1 \in [0.1,0.9]$, and add a measurement white noise $w \sim N(0,0.1)$. Again, the first 500 observations are used to estimate the models and the other 500 observations are used to validation. On Fig.~\ref{fig:res3} is shown the MSE of a three-step forecast for model 1 and the model 4. In this case, the fuzzy time series model performs better than the exact model, because the fuzzification process acts as a filter since slightly modifications in the values does not change the fuzzy relationships. The same test, for $r=3.999$ and $x_1=0.2$ is shown on TABLE~\ref{tab:res3}, where, again, the fuzzy time series model presented a smaller mean forecast error and smaller variance. 

\begin{table}[htbp]
\caption{Forecast mean squared error of the models 1 and 4 for $r=3.999$ and $x_1=0.2$ when there is a measurement noise $w \sim N(0,0.1)$ for three steps ahead}
    \centering
    \begin{tabular}{c|rr}
 & \text{Mean} & \text{Variance} \\\hline
 \text{Model 1} & \textbf{0.118824} & \textbf{0.0079419} \\
 \text{Model 4} & 0.212702 & 0.0621176 \\
    \end{tabular}
    \label{tab:res3}
\end{table}

As a fourth test, we generate 1000 observations of the logistic map with $r=4$ and $x_1=0.1$ and estimated models. After that, the models were used to forecast three steps ahead of a time series simulated with different values of $r$.

\begin{figure}[htbp]
    \centering
    \includegraphics[width=0.8\columnwidth]{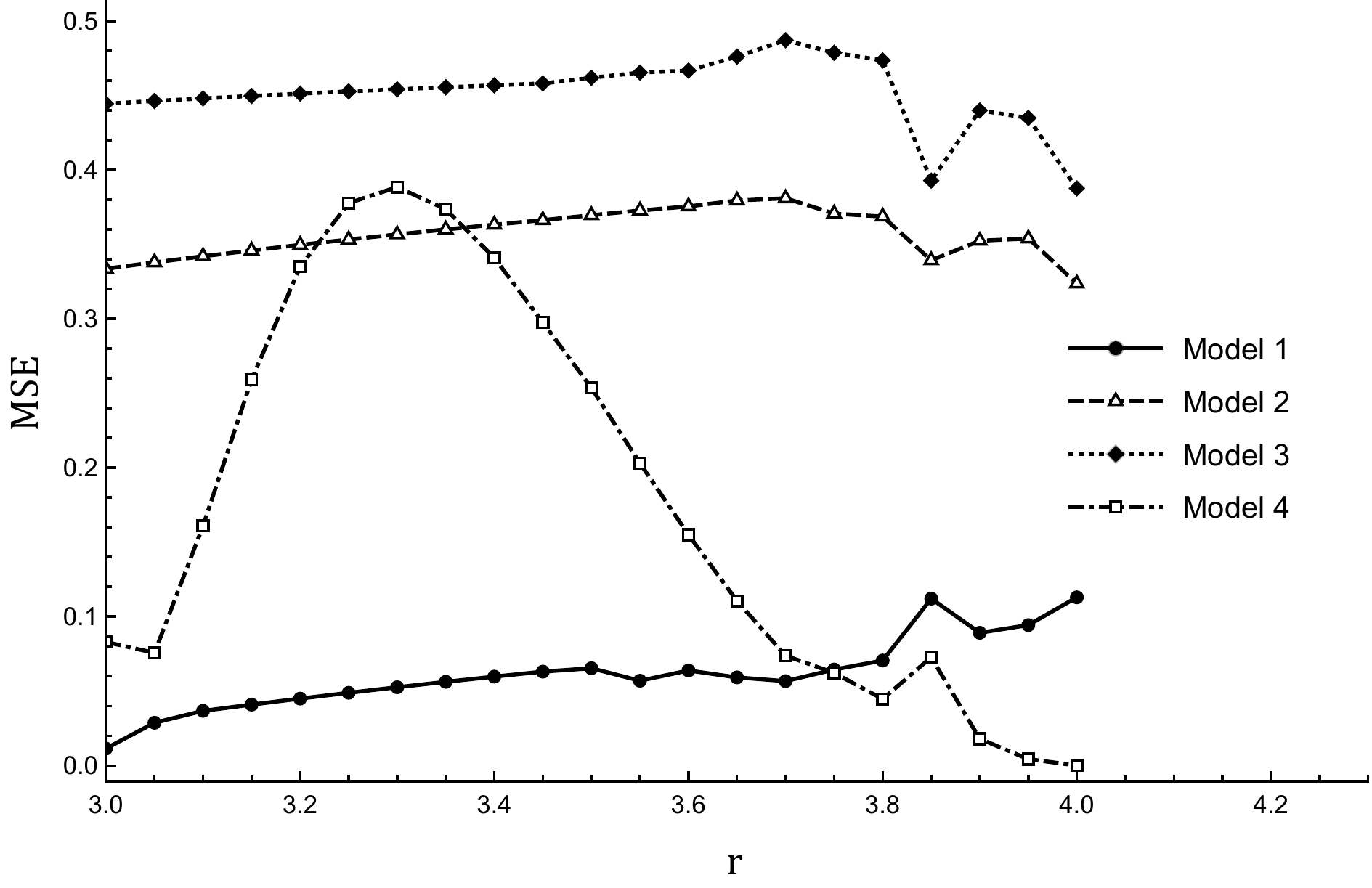}
    \caption{Forecast mean squared error of the models 1 to 4 for three steps ahead and different values $r$, with models obtained for parameter $r=4$}
    \label{fig:res4}
\end{figure}

\begin{table}[htbp]
\caption{Forecast mean squared error of the models 1 to 4 for three steps ahead and $r=3.7$ when models are obtained with data from $r=4$}
    \centering
    \begin{tabular}{c|rr}
 & \text{Mean} & \text{Variance} \\\hline
 \text{Model 1} & \textbf{0.056670} & \textbf{0.001836} \\
 \text{Model 2} & 0.378722 & 0.049714 \\
 \text{Model 3} & 0.487035 & 0.061260 \\
 \text{Model 4} & 0.072948 & 0.014470 \\
    \end{tabular}
    \label{tab:res4}
\end{table}

The results, presented on Fig.~\ref{fig:res4}, shown that the fuzzy time series model (model 1) is, also, robust to changes on the bifurcation parameter. TABLE~\ref{tab:res4} shows the same test when the test data has $r=3.7$, showing again that the model 1 presented a smaller MSE and variance, but we can observe on Fig.~\ref{fig:res4} that as the value $r$ is approaching $r=4$, the exact model (model 4) becomes improves.

\begin{figure}[htbp]
    \centering
    \subfloat[Forecast mean squared error of the models 1 (FTS) for $r=3.999$ and $x_1=0.1$ in different number of intervals]
    {
    \includegraphics[width=0.8\columnwidth]{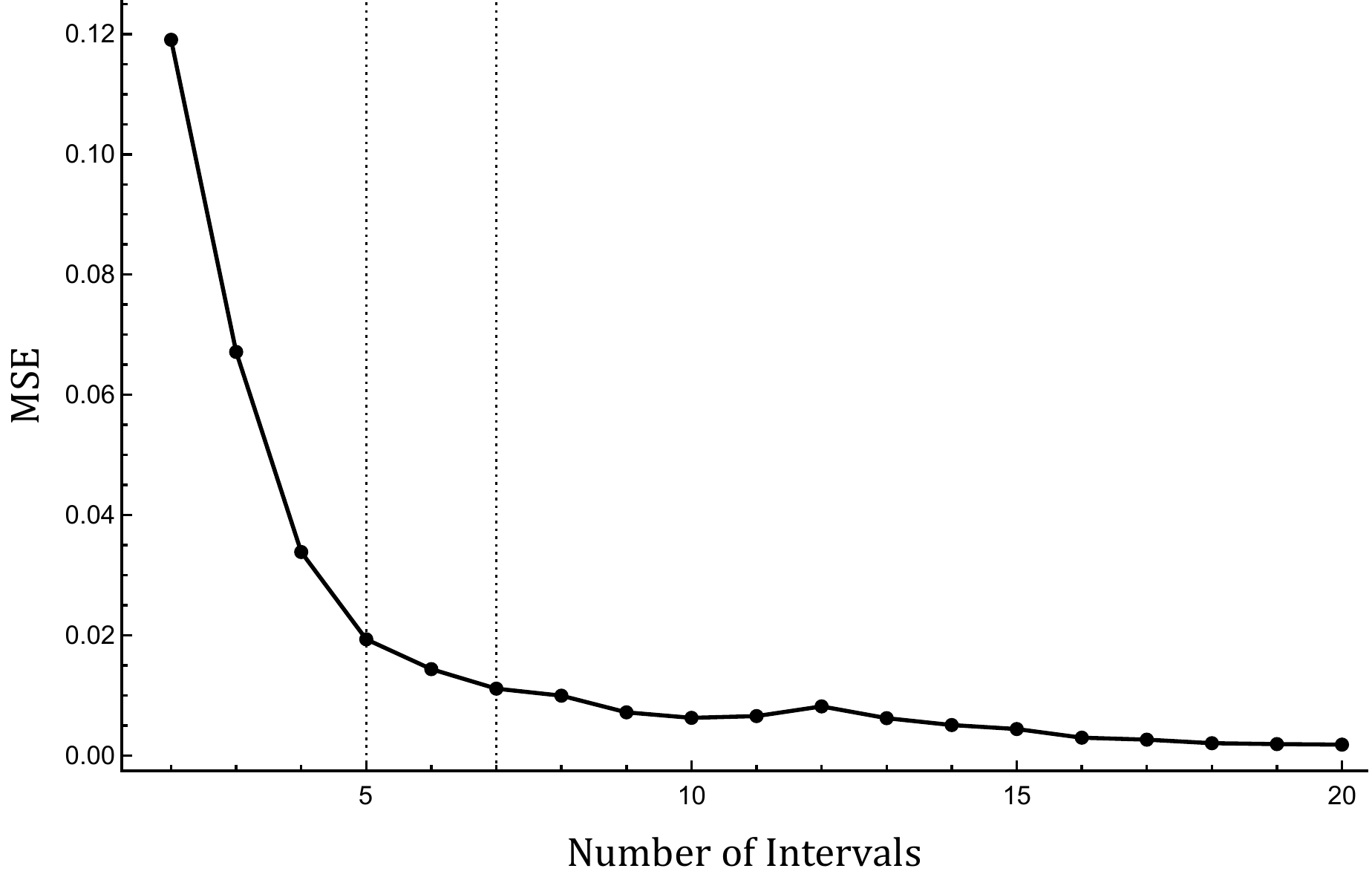}
    }
    \\
    \subfloat[Forecast mean squared error of the models 1 (FTS) for $r=3.96$ and $x_1=0.1$ in different number of intervals]
    {
    \includegraphics[width=0.8\columnwidth]{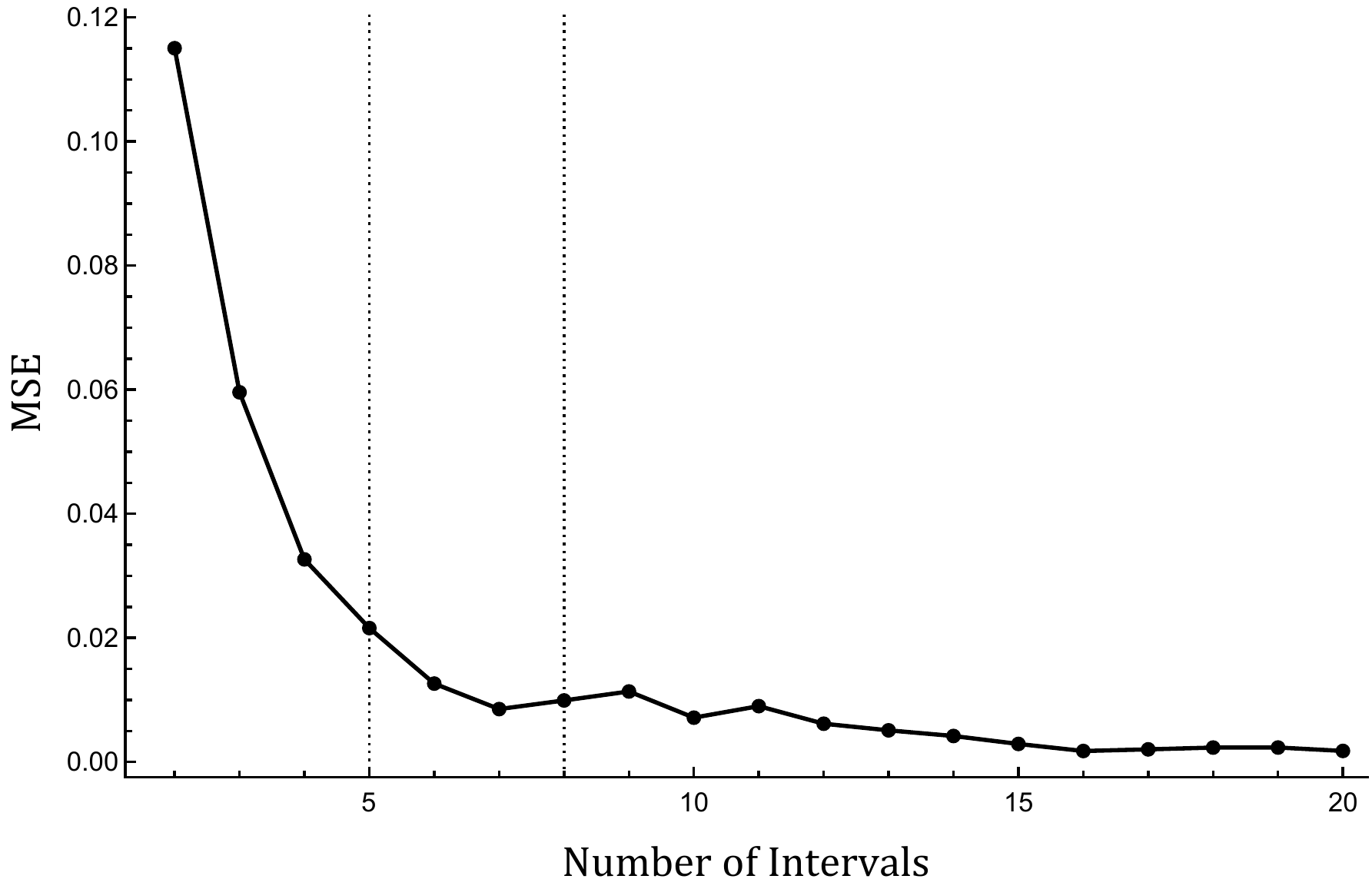}
    }
    \caption{Forecast mean squared error of the models 1 (FTS) in different number of intervals}
    \label{fig:res5}
\end{figure}

As a final test, we evaluate the Akaike Information Criterion (AIC) as an index to determine the number of intervals used to the fine the fuzzy sets. To do so, on Fig.~\ref{fig:res5} is shown the MSE of the one step ahead forecast for different number of intervals, and the number of intervals obtained by the first minimum in the AIC, and Average Method are indicated as vertical lines. On Fig.~\ref{fig:res5}(a), the AIC indicates 7 intervals and Average Method indicates 5 intervals, and increasing this number does not improve much MSE, but increases the size of the model.

On Fig.~\ref{fig:res5}(a), the AIC indicates 7 intervals and Average Method indicates 5 intervals, and increasing this number does not improve much MSE, but increases the size of the model. On Fig.~\ref{fig:res5}(b), the AIC indicates 8 intervals and Average Method indicates, again, 5 intervals. In both cases, the AIC indicates more intervals than the Average Method (in fact, in most cases the Average Method leads to 5 intervals in this problem) and generates a better model, adapting more to the specificities of each time series.  

%
\section{Conclusion}
As shown in this paper, the fuzzy time series model embeds the dynamic of the system, indicating that a first-order fuzzy time series model can forecast the logistic map, even in chaotic behavior.

The tests have shown, also, that a first-order fuzzy time series is better to forecast the logistic map than a first-order autoregressive model and a quadratic autoregressive model, and that the fuzzy times series model is less sensitive to parameter changes than the other models. In the presence of noise, the fuzzy time series model generated better forecasts than the exact model, showing its robustness.

The last test shows that the AIC is a valuable tool for the definition of the number of intervals used in the fuzzification process, despite being very conservative and, usually, picking a small number. This problem may be mitigated by manual analysis of the AIC plot, and choosing, instead of the first minimum, a greater number of intervals, which do not compromise the model size.

Fuzzy time series are efficient tools to forecast chaotic dynamical systems. As future work, we can cite the use of fuzzy time series to forecast several types of chaotic systems and improvements and an analysis of the impact of the number of intervals, the defuzzification techniques, and the use of higher-order models in the quality of the forecast

\ifCLASSOPTIONcaptionsoff
  \newpage
\fi

\bibliographystyle{IEEEtran}
\bibliography{ref}

\end{document}